

\magnification 1200
\input amstex
\documentstyle{amsppt}
\input epsf
\NoBlackBoxes
\def\fg#1{\hbox{$\pi_1(#1)$}}

\def\bd{\partial}
\def\pbd{\partial^+}
\def\nbd{\partial^-}
\def\cp{\hbox{${\Bbb C}P^2$}}
\def\ncp{\hbox{${\overline{{\Bbb C}P}^2}$}}

\def\rf{\hbox{${\Bbb R}^4$}}
\def\r{\hbox{$\Bbb R$}}
\def\Nbd{\hbox{$\Cal N$}}

\topmatter
\title{On Ribbon \rf's} \endtitle
\author{\v Zarko Bi\v zaca} \endauthor
\address{Department of Mathematics, University of Texas, Austin, TX 78712}
\endaddress
\email{bizaca\@math.utexas.edu} \endemail
\abstract  We consider ribbon \rf's, that is, smooth open 4-manifolds,
homeomorphic to \rf\ and associated to $h$-cobordisms between closed
4-manifolds.  We show that any generalized ribbon \rf\ associated to a sequence
of $h$-cobordisms between non-diffeomorphic 4-manifolds is exotic.  Notion of a
positive ribbon \rf\ is defined and we show that a ribbon \rf\ is positive if
and only if it is associated to a sequence of stably non-product
$h$-cobordisms.  In particular we show that any positive ribbon \rf\ is
associate to a subsequence of the sequence of non-product $h$-cobordisms from
[BG].\endabstract
\endtopmatter
\document
\baselineskip = 14pt \parskip = 4pt

It is well known that there are examples of pairs of homeomorphic but not
diffeomorphic simply connected closed smooth 4-manifolds, see [K].  It is also
known that any such a pair of homeomorphic, smooth, simply connected and closed
4-manifolds is $h$-cobordant, [W].  Equivalently, given such a pair of
non-diffeomorphic 4-manifolds, each one can be obtained from the other one by a
regluing of a certain open smooth 4-manifold, usually called a ``ribbon \rf'',
[DF] or [K].  A ribbon \rf\ used in a such reimbedding can be obtained from the
$h$-cobordism and is homeomorphic but not diffeomorphic to the standard
Euclidean four-space, \rf.  So, it is an example of what is usually referred as
an {\it exotic\/} \rf.  Any $h$-cobordism between a pair of smooth, possibly
diffeomorphic, 4-manifolds may be used to construct a ribbon \rf, but it is not
known whether each of them is necessarily exotic.  This paper provides a
partial answer to the question which ribbon \rf's are exotic.  We are working
under an assumption that the given $h$-cobordism can not be turned into a
product cobordism by blowing up both of its boundary components finitely many
times (see Definition 5).

Let $\ (W^5;M_0^4,M_1^4)\ $ be an $h$-cobordism between two non-diffeomorphic,
oriented, smooth, closed, simply connected 4-manifolds.  After trading handles
if necessary, we may assume that $W$ has a handlebody description with only 2-
and 3-handles:
$$W\cong (M_0\times I)\cup(\bigcup_{i=1}^k h_i^2)\cup(\bigcup_{j=1}^k h_j^3),$$
where $I$ is the unit interval, $\ I=[0,1]$, and the matrix of incidence
numbers
between 2- and 3-handles is the identity matrix.  These incidence numbers are
equal to the intersection numbers in the middle level of the cobordism, $\
M_{1/2}=\bd((M_0\times I)\cup(\bigcup_{i=1}^k h_i^2)-M_0$, between the {\it
belt\/} (or {\it dual\/}) spheres of the 2-handles and the {\it attaching}
spheres of 3-handles, see [RS].  We denote the attaching spheres by $A_i$ and
the belt spheres by $B_i, \ 1\leq i\leq k$.  (Often the belt spheres are called
``descending spheres'', and attaching spheres are called ``ascending
spheres''.)  Both the attaching spheres and the belt spheres are families of
disjointly embedded 2-spheres in $M_{1/2}$, but beside $k$ intersection points
of $\ A_i\cap B_i, \ 1\leq i\leq k$, recorded in the intersection matrix, there
may be some additional intersection points between the attaching and the belt
spheres.  These extra intersection points on any $A_i\cap B_j$ can be grouped
into pairs with opposite signs.  Note that in the absence of these extra pairs
of intersections the 2- and 3-handles in $W$ form complementary pairs of
handles that can be removed from the handlebody decomposition.  In that case
there is a product structure for $W$, that is, a diffeomorphism $W\cong
M_0\times I$.  In our situation the $h$-cobordism has no smooth product
structure, so there has to be at least one extra pair of intersections between
$A_*$ and $B_*$.

We denote by $X$ a regular neighborhood in the middle level of the union of
these spheres, $X=\Nbd (A_*\cup B_*)$.  Extra pairs of intersections result in
\fg{X} being nontrivial.  We can use Casson's construction [C] to cap the
generators of \fg{X} by Casson handles inside $M_{1/2}-\text{int}X$.  Casson's
construction may produce new pairs of intersections between $A_*$ and $B_*$,
but when considered separately, each family of spheres remains disjoint.  The
boundary components of $W$, $M_0$ and $M_1$ in our notation, can be obtained by
surgering the middle level, $M_{1/2}$.  These surgeries are performed on $A_*$
spheres to obtain $M_1$ and on $B_*$ spheres to obtain $M_0$.  Surgering $X$
produces two compacts in the boundary components, $Y_0$ and $Y_1$.  It follows
from Freedman's work [F] that $M_0$ and $M_1$ are homeomorphic and that the
cobordism $W$ is homeomorphic to the product cobordism, $M_0\times I$.  Since
we have assumed that $M_0$ and $M_1$ are not diffeomorphic, $W$ can not be
diffeomorphic to the product cobordism.  Following [DF], [FQ] or [K] we may
assume that $W$ is smoothly product over the complement of a the compact $Y_0$
in $M_0$.  Note that $\ \bd Y_0=\bd Y_1=\bd X\ $ and so $M_1$ can be obtained
from $M_0$ by replacing $Y_0$ by $Y_1$, that is $\ M_1\cong (M_0
-\text{int}Y_0)\cup_{\bd Y_0=\bd Y_1}Y_1$.  Also, $\ M_{1/2}\cong (M_0
-\text{int}Y_0)\cup_{\bd}X$.

\midinsert
\epsfxsize = 5.2in\epsfbox{y0xy1.ai}\botcaption{Figure 1}\endcaption
\endinsert
A link calculus description of an example of such $Y_0$, $X$ and $Y_1$ is
presented in Figure 1.  Dashed circles are generators of the fundamental
groups.  A compact obtained like $Y_0$ or $Y_1$ in an $h$-cobordism is known to
be diffeomorphic to a complement in the 4-ball of an embedded disc that spans a
ribbon link in the boundary of the 4-ball. We will refer to such a compact
$Y_0$ as a {\it ribbon complement associated to the $h$-cobordism $W$\/} or
simply as a {\it ribbon complement}, when the cobordism is determined by the
context.  By inverting the cobordism $W$ we obtain an $h$-cobordism $\
(\overline W;M_1, M_0)$ and $Y_1$ is a ribbon complement associated to
$\overline W$.  Meridians to the components of the bounding ribbon knot or link
generate the first homology group of a ribbon complement.  If we use these
meridians (with 0-framings) to attach the standard 2-handles to a ribbon
complement, the resulting manifold is the standard 4-ball.  If the standard
2-handles are replaced with Casson handles and the remaining boundary is
removed, then the resulting open 4-manifold is homeomorphic to the Euclidean
four-space, \rf, and is called a {\it ribbon \rf} [DF].  In the case that a
ribbon \rf\ is not diffeomorphic to \rf, we refer to it as an {\it exotic
ribbon\/} \rf.  If Casson handles are attached to ribbon complement associated
to $W$ ambiently in $M_0-\text{int}Y_0$, then we say that the resulting ribbon
\rf\ is {\it associated to the cobordism $W$}.  An example of an exotic ribbon
\rf\ associated to a non-product $h$-cobordism was explicitly described in [B].
Although only two Casson handle were involved in its construction, the number
of their kinks grow so fast with the level that the description, as I. Steward
[S] has politely phrase it, ``verges on bizzare''.  To obtain a simpler exotic
ribbon \rf, a sequence of non-product $h$-cobordisms was used in [BG] in the
following way.  Let $R_*=\text{int}Y_*\cup_{i=1}^m CH_i$ be a ribbon \rf\ built
from a ribbon complement $Y_*$. For every positive integer $n$, we denote by
$U_*^n$ the open manifold built by attaching to $Y$ only the first $n$ levels
of each of the Casson handles $CH_i$, $\ U_*^n=\text{int}Y_*\cup_{i=1}^m
(CH_i)^n$. Suppose that each $U_*^n$ is associated to an $h$-cobordism $\
(W_n;M_{n,0},M_{n,1})\ $ in the sense that $Y_*$ is associated to $W_n$ and the
$n$ level open Casson towers $(CH_i)^n$ are embedded ambiently into
$M_{n,0}-\text{int}Y_*$.  Then we say that $R_*$ is {\it associated to the
sequence of cobordism\/} $\{ W_n\} $.  Note that if $R_*$ is a ribbon \rf\
associated to an $h$-cobordism $W$ then $R_*$ is also associate to the sequence
of cobordisms $\{ W_n=W\} $.

To continue we introduce a model of a ribbon complement.  We start as before by
first constructing a compact $X$ that is a regular neighborhood of $A_*$ and
$B_*$ spheres in the middle level of an arbitrary $h$-cobordism $W$.  It is
always possible, although not in a unique way, to group the extra pairs of the
intersections so that each can be obtained by a finger move of an $A_*$ through
a $B_*$.  In other words, we may introduce finger moves on the regular
neighborhood $\ \Nbd(\coprod_k(S^2\vee S^2))\subset\sharp_k(S^2\times S^2)\ $
so that the result is diffeomorphic to $X$.  This construction produces a
distinguished set of generators for \fg{X} consisting of loops embedded into
$\bd X$.  If we retreat the fingers emanating from $A_*$ so that their tips are
only tangent to $B_*$, we call the remaining generators {\it accessory
loops\/}.  When the fingers are returned to their initial positions one of the
two intersection points of each finger is designated for accessory loops to pas
through and we adjust accessory loops accordingly.  To complete our set of
generators for \fg{X} we choose a loop for each finger that consist of an arc
on $A_*$ and an arc on $B_*$, both ending on the two intersection points on the
finger.  We call these loops the {\it Whitney loops}.  Using isotopies when
necessary, we assume that loops generating \fg{X} are disjoint outside the
extra intersections between $A_*$ and $B_*$.  After projecting $X$ along the
cobordism into $M_0$ or $M_1$, only the interior of $X$ has been replaced and
the accessory and the Whitney loops are in $M_0$ or $M_1$, respectively.  It is
easy to describe this projection in the terms of link calculus: to surger $X$
into $Y_0$ 0-framings of the link components representing the family $B_*$ are
replaced by dots, namely, 2-handles are replaced by 1-handles or by scooped out
2-handles.  Note that \fg{X} is a subgroup of \fg{Y_0}, the latter also
includes generators that are meridians to the dotted circles that used to
represent $B_*$ in $X$.  We will continue to call ``accessory'' and ``
Whitney'' loops the generators for \fg{Y_0} induced by the surgery.

\midinsert
\epsfxsize = 5in\epsfbox{finger.ai}\botcaption{Figure 2}\endcaption
\endinsert

Figure 2 is a link calculus picture of a finger move.  Note that the components
of $B_*$ are already surgered into dotted circles, the picture of $X$ can be
obtained by replacing dots on these components by 0-framings.  We may build our
model of a ribbon complement by adding such finger moves to a collection of
Hopf links with a 0-framed and a dotted component in each (link calculus
picture of complementary pairs of 1- and 2-handles), but often we will end up
with more ``accessory'' dotted circles than needed.  The extra dotted circles
may me slid of their parallels and off our picture where they are removed.
Alternatively, the ``accessory'' dotted circle from Figure 2 is added only when
the finger closes a loop and introduces a new accessory generator of \fg{X}.

\midinsert
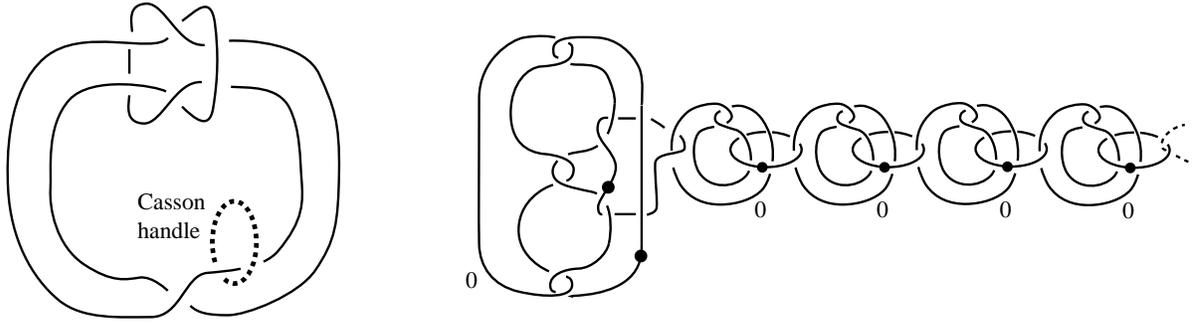
\epsfxsize = 5.2in\epsfbox{bg1.ai}\botcaption{Figure 3} $R_0$ -- an example of
a generalized ribbon \rf.\endcaption
\endinsert
\remark{Remark} The construction of ribbon complements described above (and in
[K] and [DF]) involves a specific family of ribbon links.  The above mentioned
exotic \rf\ from [B] is an example of an ribbon \rf, but the simplest known
exotic \rf\ (introduced in [BG]) is not a ``ribbon \rf'' although it has a
ribbon knot (Figure 3) associated to it: to the meridian of the ribbon
complement of the disc bounding ribbon knot from the left part of Figure 3 a
single Casson handle is attached. The attached Casson handle has a single kink
at each level and all its kinks are positive.  This exotic \rf, which we denote
by $R_0$, is built from the same ribbon complement as the example from [B], but
instead capping both the accessory and the Whitney loops by Casson handles,
there are only one Casson handle and a standard 2-handle involved.  We will
consider a slightly more general situation and so we say that a contractible
open smooth 4-manifold built from a ribbon complement by capping the accessory
and Whitney loops with any combination of 2- and Casson handles is a {\it
generalized ribbon \/} \rf.  The notion of such a generalized ribbon \rf\ being
associated to a sequence of $h$-cobordisms can be defined exactly as
before. \endremark

Two questions arise naturally.  The first one is whether a ribbon \rf\
associated to an $h$-cobordism between non-diffeomorphic 4-manifolds (or to a
sequence of such cobordisms) is necessary exotic.  Conversely: which
combinations of ribbon complements and attached Casson handles produce exotic
ribbon \rf's.  Answer to the first question was known to be positive in the
case of a ribbon \rf\ associated to a single $h$-cobordism between
non-diffeomorphic 4-manifolds, [K, pages 98 -- 101] .  This is also true in a
slightly more general situation.

\proclaim{Theorem 4}{\rm (Compare with Theorem 3 in [K, page 98].)}  A
(generalized) ribbon \rf\ associated to a sequence of $h$-cobordisms between
non-diffeomorphic 4-manifolds is exotic.\endproclaim

\demo{Proof} Let $R=Y\cup\big (\bd Y\times (0,\infty)\big)\cup_{i=1}^m
CH_i\cup_{j=1}^n H^2_j\ $ be a generalized ribbon \rf\ associated to a sequence
of $h$-cobordism, $\{ (W_n;M_{n,0},M_{n,1})\}$, where $M_{n,0}$ and $M_{n,1}$
are not diffeomorphic and where $H_j^2$ denotes an open 2-handle, that is $\
(H_j^2,\bd H_j^2)\cong (D^2\times\r,S^1\times\r)$.  Assume that $R$ is
diffeomorphic to the standard \rf.  Then, since the ribbon complement $Y$ is a
compact subset of $R$, there is a smooth 4-ball $B_0$ embedded in $R$ that
contains $Y$ in its interior.  The ball $B_0$ being compact is contained in
$U^k$, for some $k\geq 1$.  $W_k$ is a product over $\ M_{k,0}-Y\ $ and we use
this product structure to lift $\bd B_0$ into a smooth 3-sphere $S_1$ in
$M_{k,1}$.  An argument from [K] shows that $S_1$ has to bound a standard
4-ball in $M_{k,1}$: briefly, by embedding $k$ level Casson towers of $U^k$
into the standard 2-handles we construct a smooth embedding of $U^k$ into the
standard 4-ball which we consider to be in the standard 4-sphere, $S^4$.  The
piece of the cobordism $W_k$ over $U^k$ can be transplanted into the product
cobordism $S^4\times I$.  Recall that the complement of a smooth 4-ball in
4-sphere is also a 4-ball.  We lift the complement of the 4-ball $B_0$ in $\
S^4\times\{ 0\}\ $ to the top of the cobordism, $\ S^4\times \{ 1\}$. Since
this cobordism is a product over the complement of $Y$, the complement of the
lifted 4-ball is a standard 4-ball, bounded by $S_1$, and the cobordisms over
the complements of int$B_0$ in $S^4$ and $M_{k,0}$ are product.  So, the
product structure over $\ M_{k,0}-B_0\ $ can be extended over $B_0$,
contradicting our assumption that $M_{k,0}$ and $M_{k,1}$ were not
diffeomorphic.\qed\enddemo

One might expect that the second question has an equally simple answer, that
all possible generalized ribbon \rf's are exotic.  This is not the case: the
simplest ribbon complement is diffeomorphic to $\ S^1\times D^3$ and the ribbon
link involved is the unknot.  If any Casson handle is attached over the
meridian and the boundary is removed, the resulting manifold is the standard
\rf\ [F, page 381].  However, it is easy to describe this manifold as a
generalized ribbon \rf, for example, replace one of the two Casson handles of
$R_1$ in Figure 5 by a standard open 2-handle.

\definition{Definition 1} If $L$ is an accessory loop of a ribbon complement
$Y_0$ than the set of Whitney loops that intersect $L$ is called the {\it
Whitney set of\/} $L$.
\enddefinition

\definition{Definition 2} If a signed tree associated to a Casson handle has a
positive branch, than the Casson handle is {\it positive\/}.  If there are more
positive then negative edges emanating from every vertex of the associated
tree, then the Casson handle is {\it strictly positive\/}.
\enddefinition

\definition{Definition 3} Let $Y$ be a ribbon complement and $R$ a ribbon \rf\
obtained from $Y$ by adding Casson handles.  Suppose that there is an accessory
loop such that every loop of its Whitney set is capped by a positive Casson
handle and, in the case that its Whitney set contains only one loop, then the
accessory loop itself is capped by a positive Casson handle. In the case that
there are more then one loop in the Whitney set we require that this accessory
loop coincides with at most one finger emanating from any $A_*$ spheres.  Then
we say that $R$ is a {\it positive ribbon \/} \rf.
\enddefinition

It is not known whether ribbon \rf's that are not positive are exotic or not.
However, ribbon \rf's that are not positive can not be associated to a
sequence of non-product $h$-cobordisms that satisfy the following additional
assumption.

\definition{Definition 5} Let $(W^5;M_0^4,M_1^4)\ $ be an  $h$-cobordism
between two oriented, smooth, closed, simply connected 4-manifolds.  We say
that $W$ is {\it stably non-product\/} if $M_0\sharp n(\ncp)$ and $M_1\sharp
n(\ncp)$ are not diffeomorphic for any nonnegative integer $n$.
\enddefinition

\remark{Remark} It is not known to the author whether there exists a pair of
simply connected closed smooth 4-manifolds $M_0$ and $M_1$ that are
homeomorphic and non-diffeomorphic and such that $M_0\sharp n(\ncp)$ and
$M_1\sharp n(\ncp)$ are diffeomorphic for some $\ n\geq 1$.
\endremark

\proclaim{Theorem 6} A ribbon \rf\ associated to a sequence of stably
non-product $h$-cobordisms is positive.\endproclaim

We will prove that a ribbon \rf\ that is not positive can not be associated to
a sequence of stably non-product $h$-cobordisms $W_n$ by showing that at least
one $W_n$ can be turned into a product cobordism by blowing up its end
sufficiently many times.  A process of removal of double points is described in
[Ku] and a short outline of that method is given next.

Suppose that $(\Delta,\bd\Delta)$ is an immersed disc in a 4-manifold $(N,\bd
N)$ with a single double point in the interior of $\Delta$.  Furthermore,
suppose that this double point is negative.  After blowing up $N$ we can
replace $\Delta$ by an embedded disc in $(N\sharp\ncp,\bd (N\sharp\ncp))$ that
spans the same loop in $\bd N\sharp\ncp=\bd N$: Let $E$ and $E'$ represent two
``exceptional curves'', i.e., copies of ${\Bbb C}P^1$ in general position and
embedded in the added \ncp. If $E$ and $E'$ are equipped with opposite
orientations then they intersect in a single point that has the positive sign.
Choose a small ball centered at the intersection between $E$ and $E'$, the
intersection between the boundary of the small ball and $E$ and $E'$ will form
a Hopf link.  Similarly choose a small ball in the interior of $N$ that is
centered at the double point of $\Delta$.  The intersection between $\Delta$
and the boundary of the ball is again a Hopf link.  The centers of these two
balls can be connected by a path that avoids $\Delta$, $E$ and $E'$.  Remove
the intersection between the interiors of the balls and $\Delta$, $E$ and $E'$.
Now the two Hopf links in the boundaries of the balls are connected by two
pipes that follow the chosen path.  The resulting disc represents the same
second homology class as $\Delta$ and has the same boundary, but it is embedded
into $(N\sharp\ncp,\bd (N\sharp\ncp))$.

Notice that this procedure can prune all the branches of a tree associated to
a Casson handle that have a negative kink.  Also, for every non-positive
Casson handle there is a natural number $k$ so that every brunch of the tree
associated to the Casson handle has a negative kink on the first $k$ levels.
(Recall that a tree associated to a Casson handle has finitely many edges
coming from every vertex.)  Consequently, if we perform sufficiently many
blow-ups of an ambient 4-manifold, we may replace any Casson handle that is
not positive with an embedded standard 2-handle.

\demo{Proof of Theorem 6} Suppose that $R$ is a non-positive ribbon \rf\
associated to a sequence of $h$-cobordisms $W_n$.  Choose $k$ large enough
such that every branch of non-positive Casson handles contains a negative kink
on the first $k$ levels.  We will work in $W_k$ where we have embedded first
$k$ levels of the Casson handles from $R$.  After blowing up $W_k$
sufficiently many times we may replace all non-positive Casson towers in $\
\tilde W_k :=W_k\cup(\sharp\ncp\text{'s}\times I)\ $ with standard 2-handles
so we may assume that all the remaining Casson towers have only positive kinks.
We can use the embedded 2-handles to perform Whitney tricks and cancel pairs of
2- and 3-handles from the handlebody decomposition of $\tilde W_k$ whenever
possible.  If none of the Casson handles of $R$ is positive, this procedure
removes all the 2- and 3-handle pairs and $\tilde W_k$ has a product structure.
In particular, $W_k$ is not stably non-product.
\midinsert
\epsfxsize = 4in\epsfbox{norman.ai}\botcaption{Figure 4}\endcaption
\endinsert

If there are no accessory loops we can use Norman trick from Figure 4 to remove
the extra pairs of intersections.  Working in the middle level we start from
such an extra pair, say between $A_i$ and $B_j$.  Note that $B_j$ and $A_j$
have a single intersection point since otherwise there would be an accessory
loop on them.  Each of the two intersections is removed by Norman trick (see
[FQ]) by removing a disc from $A_i$ centered at the intersection point and a
disc from a copy and $A_j$ that is centered at an intersection point between
$A_j$ and $B_j$.  Then the boundaries of these two removed discs are meridians
to $B_j$ and are connected by a tube, Figure 4.  The resulting new $A_i$
sphere, denoted by $A'_i$, intersect each $B_*$ that $A_j$ did, so there are
four intersections between $A'_i$ and $B_k$ in our example from Figure 4.
Repeating this process produces a cascade of fingers, each piercing a $B_t$
such that $A_t$ has no fingers. The application of Norman trick will add two
copies of $A_t$ to each finger that ends on $B_t$ therefore removing the pair
of intersections we have started with.  Now all Whitney discs are removed and
again, $\tilde W_k$ has a product structure.

If $Y$ does contain accessory loops, since we have assumed that $R$ is not
positive, each accessory loop ventures over more then one finger emanating from
a single $A_*$ or the Whitney set for the accessory loop contains a loop capped
by a non-positive Casson handles.  In the later case, after the blow-ups this
Whitney loop is capped by the standard two handle and the finger containing it
is removed, breaking the given accessory loop.  So now we may assume that the
only accessory loops remaining are those that contain more then one finger
emanating from a single $A_*$.  Starting with two such fingers and pushing them
over other spheres from $A_*$ we produce cascades of fingers.  The accessory
loop is closed when both cascades of fingers intersect the same $B_*$.  We
consider such a loop, emanating from, say, $A_i$ and ending on $B_j$.  Fingers
that may start from $A_j$ also can be grouped in pairs, each ending on a same
$B_*$.  Following these pairs we can not close a loop by having a pair of
fingers ending on $B_i$, otherwise we could select a member from each pair and
obtain an accessory loop that passes over at most one finger emanating from any
$A_*$.  In that case the non-positiveness would apply that at least one of the
associated Whitney loops can be capped by the standard 2-handle.  Consequently
we may assume that our pairs of fingers emanating from $A_i$ ends on a $B_j$
such that $A_j$ has no fingers and no extra intersection points. Using the
Norman trick as before the pairs of intersections on fingers are removed.  So
all the extra pairs of intersections can be removed and again, $\tilde W_k$ has
a product structure.\qed\enddemo

The converse to Theorem 6 is also true, every positive ribbon \rf\ is
associated to a sequence of stably non-product $h$-cobordisms.

\proclaim{Theorem 7} Every positive ribbon \rf\ can be associated to a
subsequence of the sequence $\ \{W_m \}^{\infty}_{m=2}\ $ of stably non-product
$h$-cobordisms constructed in [BG].  \endproclaim

Each of these $h$-cobordisms from [BG] we denote here by
$(W_m;M_{m,0},M_{m,1})$, where $M_{m,1}\cong E(m)\sharp k(\ncp)$ and $M_{m,0}$
decomposes as a connected sum of \cp's and \ncp's.  For simplicity we have not
included ``$k$'' or $\tilde W_m$ in our notation and $E(m)$ denotes the minimal
elliptic surface with no multiple fibers and of the Euler characteristic $12m$.

\midinsert
\epsfxsize = 4.5in\epsfbox{universal.ai}\botcaption{Figure 5}\endcaption
\endinsert
\demo{Proof} First we will show that for each positive ribbon \rf,
$R=\text{int}Y\cup_i CH_i$ and for each natural number $k$ we can embed $\
U^k=\text{int}Y\cup_i (CH_i)^k$ in some $M_{m,1}$, when $m$ is large enough and
$M_{m,1}$ contains sufficiently many copies of \ncp.  Each of this embeddings
factors through an embedding into a compact obtained from the closure of
$U_0^k$ (that is, the first $k$ levels of $R_0$ from Figure 3) by adding extra
2-handles and parallel copies of the Casson tower, shown in Figure 7.  Then we
show that the $h$-cobordism obtained by regluing the embedded ribbon complement
$Y$ is diffeomorphic to $W_m$.

\midinsert
\epsfxsize = 5in\epsfbox{iso1.ai}\botcaption{Figure 6}\endcaption
\endinsert
According to Definition 3 each positive ribbon \rf\ contains an accessory loop
whose Whitney set is capped by positive Casson handles.  We start our
construction by embedding all the other Casson handles into the standard
2-handle and each positive Casson handle capping a loop from the fixed Whitney
set is embedded into the $CH^+$, the positive Casson handle with one kink per
level.  Now each positive ribbon \rf\ is embedded into a ribbon \rf\ that has a
single accessory loop, similar to one of those in Figure 5.  In this figure
there is a sequence of ribbon \rf's, $R_n,\ n\geq 1$, each but the first one is
built by attaching $n$ copies of the Casson handle $CH^+$ to a ribbon
complement which we denote by $Y_n$.  Note that in $R_1$ the accessory loop is
also capped by $CH^+$, and for $n\geq 2$, the accessory loop is capped by the
standard 2-handle and it's dotted circle disappears (compare with Figure 2).
We denote by $R_n'$ the resulting ribbon \rf\ in which we have embedded $R$.
The index $n$ in ``$R_n'$'' or ``$R_n$'' is equal to the number of pairs of $A$
and $B$ spheres that form the underlying middle level compact, $X_n$.

The middle level compact $X_n$, and therefore the ribbon complement $Y_n$, are
uniquely defined up to isotopy by listing the geometric and algebraic numbers
of intersections between $A_*$ and $B_*$ spheres.  Furthermore we can isotope
one link calculus picture of $Y_n$ into another by sliding the 2-handles and
dotted circles whose meridians are Whitney loops over the dotted circle
corresponding to the accessory loop.  Equivalently, the possible differences
between link calculus pictures may occur as different choices of clasps,
positions of dotted circles whose meridians are Whitney loops and twists of
parallel strands in 2-handles.  In Figure 6 it was shown how to deal separately
with each of this differences.  Since in any stage we are allowed to blow-up
the ambient manifold finitely many times we can always introduce positive
twists by attaching a $-1$-framed 2-handle and sliding other 2-handles off it,
Figure 6.  So we may assure that all clasp between handles corresponding to
$A_*$ and $B_*$ spheres and positions of the dotted circles corresponding to
Whitney loops in a link calculus picture of $R_n'$ are exactly the same as in
the picture of $R_n$ in Figure 5.  The only possible difference remaining is in
accessory loops.  Therefore, our construction produces an embedding of a
positive ribbon \rf\ into a possibly blown-up ribbon \rf\ that we have denoted
by $R'_n$ and that is built by attaching copies of the Casson handle $CH^+$
onto Whitney circles of $Y_n$, and by dealing with the accessory loop in the
same fashion as in $R_n$.  If we fix as generators for \fg{Y_n} the Whitney
circles $\ w_1,\dots,w_n\ $ and the accessory loop $a$ from Figure 5, then in
general $a'$, the accessory loop of $R'_n$, is a word in these generators
involving $w_i$'s.  As before $U^k_n$ will denote $Y_n \cup_{w_*}
n(CH^+)^k)\cup_a h^2$, where $(CH^+)^k$ is the Casson tower equal to the
closure of the first $k$ levels of the Casson handle $CH^+$.  Copies of
$(CH^+)^k$ are attached over the Whitney loops $w_*$ and the 2-handle $h^2$ is
attached over the accessory loop $a$.  Similarly, we define $(U')^k_n$ by
attaching the 2-handle over $a'$ instead of $a$ .

\midinsert
\epsfxsize = 3in\epsfbox{container.ai}\botcaption{Figure 7}\endcaption
\endinsert
We claim that we can extend the embeddings from [BG] of $U^k_0$ into $M_{m,1}$,
where $m$ is large enough, to an embedding of the handlebody $C_1$ from Figure
7.  The embeddings we are extending were described in [B], Figures 19 -- 60,
and in [BG], Figures 40 -- 81.  The modifications we have to add is to have an
arbitrary number of Casson $k$ level towers (instead of only one in $U^k_0$)
and arbitrary numbers of $-1$-framed 2-handles linked with the two dotted
circles in $C_1$, Figure 7.  To embed $-1$-framed 2-handles linked with the
larger dotted circle we follow Figures 40 -- 43 and 49 -- 54 in [BG].  In
Figure 40 from [BG] the meridian of the larger dotted circle from our Figure 7
corresponds to the circle denoted by $\beta$.  The meridian of the smaller
dotted circle in Figure 7 we can follow in Figures 39 -- 44 in [B], but the
difference in our case that the actual pictures are the mirror images of those
in [B] and the largest 0-framed two handle and the dotted circle have to switch
their roles.  So we read Figure 44 in [B] that each meridian of the smaller
dotted circle form our Figure 7 is isotoped to a pair of meridians of the
larger dotted circle from Figure 7.  Each pair has one unlinked 0-framed
component and the other one is $-1$-framed and linked with all the second
components of the other pairs.  After passing these meridians into the other
part of the manifold we have to take their mirror images so the framings in
Figure 44 from [B] are now correct as drown.  Figure 60 from [B] shows how to
deal with (now $+1$-framed) linked second components of the pairs.  Note that
in figures from [BG] these pairs of meridians also are isotopic to $\beta$.  In
all cases we are left to cap $0$-framed isotopes of the circle $\beta$ which we
can do by either Casson towers of arbitrary levels (Figure 59 in [B] or Figure
81 in [BG]) or we can slide them over the linked $-1$-framed 2-handles and
produce embeddings of $-1$-framed 2-handles in our Figure 7.  Each of these
processes uses $-1$-framed 2-handles isotop to $\beta$ (Figure 81 from [BG])
end to procure them in a sufficient quantity we need to choose $m$ large
enough.

\midinsert
\epsfxsize = 5in\epsfbox{embedding.ai}\botcaption{Figure 8}\endcaption
\endinsert
Next we construct an embedding of $U^k_n$, $n\geq 1$, into the handlebody $C_1$
from Figure 7.  An embedding of $U^k_1$ was described in [BG], Figure 47,
namely a 0-framed 2-handle is added to connect the accessory and Whitney dotted
circle, the result is $U^k_0$, but with two parallel Casson $k$ level towers so
it embeds in $C_1$.  To embed $U^k_n$, $n\geq 2$, into $C_1$ we connect the
dotted circles corresponding to $B$ spheres by $n-1$ 2-handles.  Figure 8
depicts this process in the case of $U^k_2$ and $U^k_3$.

\midinsert
\epsfxsize = 5in\epsfbox{inside.ai}\botcaption{Figure 9}\endcaption
\endinsert
Figure 9 shows how to embed each $U^k_n,\ n\geq 1$, in $C_1$.  We start with
$C_1$ and add complementary pairs of 2- and 3-handles (which corresponds in
link calculus to adding unlinked unknots with framings 0) and then we slide the
2-handles over the 2-handle of $C_1$.  Next we perform isotopies to separate
these parallel copies of this 2-handle.  Figure 9 shows how to complete an
embedding of $U^k_2$ into $C_1$ and all the other $U^k_n$'s are embedded in the
same fashion.

\midinsert
\epsfxsize = 5in\epsfbox{iso2.ai}\botcaption{Figure 10}\endcaption
\endinsert

Next we show how to modify these embeddings to embed $(U')^k_n$ into $C_1$.
Recall that the difference between $(U')^k_n$ and $U^k_n$ is only in the
2-handle capping the accessory loop.  In the case of $(U')^k_n$ the attaching
circle of this 2-handle can also link other handles, but we can isotop this
circle to be a word in our fixed Whitney circles, $w_1, w_2,\dots,w_n$.  Figure
10 shows how to unlink the accessory loop from the dotted circles corresponding
to Whitney loops.  For each piece of the accessory loop that links once a
Whitney dotted circle we add a $-1$-framed 2-handle, as shown in Figure 10.
Then we can use the added handle to slide the accessory loop off the dotted
circle.  The result of such a slide increases the framing of the accessory loop
by $+1$.  After we slid the accessory loop off all dotted circles the resulting
new accessory loop is linked with only one dotted circle, as in the case of
$U^k_n$, but its framing will be in general positive.  We can slide in Figure
10 the dotted circle and the $N$-framed accessory circle linked to it such that
the dotted circle ends up linked with the visible 0-framed 2-handle. Then we
can slid each of the two strands of the 0-framed handle over the $N$-framed
handle and finally, we slid this (still 0-framed) 2-handle $2N$ times over
dotted circle to unlink it from the $N$-framed handle.  The result of this
handle slides will render the canceling pair with the dotted circle and the
$N$-framed 2-handle unlinked from anything else and therefore it can be removed
from the picture.  The other induced change is that the two strands of the
0-framed 2-handle have obtained $N$ positive twists.  The framing of the
accessory loop can always be increased by any positive amount: attach a
$-1$-framed 2-handle as in Figure 10 and slide off it the 0-framed handle.
This process blows up once the ambient manifold and introduces an extra
positive twist.  By using such a blow-up if necessary and reversing the handle
slide, we assume that $N$, the framing of the accessory circle, is an even
positive integer.

\midinsert
\epsfxsize = 5in\epsfbox{trick.ai}\botcaption{Figure 11}\endcaption
\endinsert
Our present link calculus pictures of $(U')^k_n$ and $R'_n$ differ from $U^k_n$
and $R_n$ in that accessory loop is capped by an unknoted 2-handle that has
framing $N$ instead of 0, where $N$ is an even positive integer or,
equivalently, one of the 0-framed 2-handles has $N$ positive twists, see Figure
10.  Figure 11 shows how to embed $(U')^k_n$ into $C_1$ by accommodating each
pair of twists.  Adding a pair of complementary 1- and 2-handles such that the
2-handle is 0-framed and linked twice with the dotted circle is equivalent of
having two positive twists and so we add $N/2$ such pairs.  The attaching of
$-1$-framed 2-handles as in Figure 11 replaces each of $N/2$ pairs of twists by
a $-1$-framed 2-handle that is meridian to the larger dotted circle in the
picture of $C_1$ and so we have an embedding of $(U')^k_n$ into $M_{m,1}$.

\midinsert
\epsfxsize = 5in\epsfbox{embed2.ai}\botcaption{Figure 12}\endcaption
\endinsert
Our next task is to obtain an $h$-cobordism from each embedding of $(U')^k_n$
into $M_{m,1}$.  The other boundary component, $M_{m,0}$, is obtained by a
reimbedding of $\ (Y')_n\subset (U')^k_n \subset M_{m,1}\ $ that switches the
roles of 0-framed 2-handles and appropriate dotted circles in Figures 8 and 9.
In particular, Figure 9 is replaced by Figure 12.  The 0-framed Hopf links
added in Figure 12 were complementary pairs of 1- and 2-handles in Figure 9.

We claim that the obtained $h$-cobordisms are $(W_m;M_{m,0},M_{m,1})$, that is,
they are in the same sequence of $h$-cobordisms obtained by the reimbeddings of
$\ Y_0 = Y_1\subset U^k_1 \subset M_{m,1}$, [BG], but with possibly different
$m$ and the number of blow-ups corresponding to a given number of levels, $k$.
The boundary components of the cobordisms in [BG] were constructed by regluing
of the Mazur rational ball that is visible in our figures as a subhandlebody of
$C_1$, Figure 7.  The two embeddings of the Mazur ball were using its dual
handlebody decomposition and we will do the same here with embeddings of $C_1$.
We will explicitly show the embedding of $Y_2$ into $M_{m,0}$ and from the
construction it will be clear how to obtain the embeddings of $Y_n, n\geq 3$.

\midinsert
\epsfxsize = 5.2in\epsfbox{dual1.ai}\botcaption{Figure 13}\endcaption
\endinsert
The top part of Figure 13 recapitulates the embedding of $Y_2$ into $C_1$.
There is a 3-handle attached over unlinked 2-handle, normally not visible in a
link calculus picture.  Below is a link calculus picture of the dual handlebody
decomposition.  We will follow a convention from [BG] and now we present its
shout outline.  To obtain a link calculus picture of a dual decomposition from
a given link calculus picture one may start by drawing the mirror image of the
given link calculus picture.  Then the dotted circles (that is, 1-handles or,
equivalently, scooped out 2-handles) obtain $(0)$-framings and the signs of all
others framings are changed and enclosed by parentheses, see [BG].  Next, to
each link component that was an attaching circle of a 2-handle in the original
link calculus picture one attaches a 0-framed 2-handles over the meridian of
the component, see [K].  (Recall that in a dual decomposition of a
4-dimensional handlebody 1-handles become 3-handles and vice versa and the
0-handle becomes a 4-handle.)  Such a link calculus picture contains components
marked with a dot (1-handles), components with an integer framing (2-handles)
and components whose framings are integers enclosed by parentheses.  A
handlebody described by a such link has two boundary components: the
``$\nbd$-component'' is obtained by performing (a 3-dimensional) surgery only
on components in parentheses, and the other boundary component, the
``$\pbd$-component'', is the result of the surgery of all the components of the
link.

To facilitate the description we have labeled all components of the dual part
of Figure 13 by capital letters, A -- I.  Now we describe a diffeomorphism
between the two pictures of the dual decomposition.  The handlebody to the left
also contains a 4-handle and three 3-handles, the duals of the 1-handles in the
original decomposition, namely they have to be attached over the components D,
E and G.  Also there is a 1-handle that has to be added to the $\nbd$-boundary
component that is not visible in this picture.  First we slide D and E over I
and off G.  The new D' and E' are now linked with A, and A can be slid off I
over G.  Components I and G are unlinked and can be removed from the picture.
Note that now A occupies the place previously occupied by G. A is then slid of
D and E over B and C, respectively.  The resulting 2-handle, now denoted by A'
is unlinked from the rest of the components, but the 3-handle that used to be
attached over G is now attacher over A' and together they form a complementary
pair of 2-and 3-handles.  Next, we slide D' over E' and the result is visible
in the lower right corner of Figure 13.  Now the 1-handle is visible; it
coincides with D' and together with the 2-handle B it forms a complementary
pair.  Now we have two complementary pairs of handles that we remove from the
picture.  By adding where appropriate Casson towers and $-1$-framed 2-handles
we obtain $C_1$.

\midinsert
\epsfxsize = 5.2in\epsfbox{dual0.ai}\botcaption{Figure 14}\endcaption
\endinsert
We proceed with a description of an embedding of $Y_2$ into $M_{m,0}$.  The
top of Figure 14 reproduces from Figure 12 the link calculus picture of $Y_2$
with two 0-framed 2-handles added.  Below is a link calculus picture of its
dual decomposition.  As before we have labeled all the components of this link.
Again we have an invisible 1-handle, a 4-handle and three 3-handles attached
over D, E and F.  First we slide B over E and then twice over C to unlink it
from H.  The resulting component, B', is unlinked from the rest of the
components, but it coincides with the 3-handle attached over E.  Similarly as
above, we remove the components G and I and the resulting link calculus picture
is visible in Figure 14 and the next picture in that figure is the result of an
ambient isotopy of the link.  Then we slide E' over D' and the resulting
component, E'', is where we add (from ``below'', to the $\nbd$-component of the
boundary) the missing 1-handle.  The 3-handle that was incident only with E is
now incident with both E'' and B' and 3-handle originally attached to D is now
incident with D' and E''.

\midinsert
\epsfxsize = 4in\epsfbox{dual02.ai}\botcaption{Figure 15}\endcaption
\endinsert
The $\nbd$-component of the boundary of our handlebody is the same as in the
case of embedding into $M_{m,1)}$ and, after adding two $-1$-framed 2-handle we
have a handlebody from Figures 44 and 45 in [BG].  The result of this addition
is in Figure 15.  Furthermore we decompose this manifold as a union of three
pieces stack over each other and glued over appropriate boundary components.
The bottom piece in Figure 15 is the dual picture of adding two $-1$-framed
2-handles.  Its $\pbd$ boundary component is obtained by surgeries on D', F and
H.  Then, we add C, a 0-framed 2-handle.  Now, the $\pbd$ boundary component of
this middle piece is $S^3$ and the difference between the handlebody we are
considering and the one used in embedding of $Y_2$ into $M_{m,1}$ from [BG] is
in the pair of 1-handle and a 0-framed 2-handle added in the piece on the top
and glued by a diffeomorphism of the standard 3-sphere.  Since by changing a
gluing diffeomorphism of the standard 3-sphere we can not change the smooth
structure of the resulting 4-dimensional manifold, we have only to consider the
diffeomorphism type of the piece on the top.  It is easy to see that this
piece, together with 3-handles and the 4-handle is diffeomorphic to the
standard 4-ball.  Namely, by sliding D' over $(0)$-framed C we can unlink all
the $(0)$-framed components and then E'' and A' form a complementary pair of 1-
and 2-handles.  Now the 3-handles are attached over the resulting $(0)$-framed
components, and on top of them we have to attach the 4-handle.  Therefore, the
handlebody from Figure 15, together with invisible 3- and 4-handles is the same
as the manifold from Figures 44 and 45 in [BG] and the argument there shows how
to decompose the obtained boundary component of the $h$-cobordism, $M_{m,0}$,
into a connected sum of \cp's and \ncp's.  To generalize to embeddings of
$Y_n,\ n\geq 3$ into $M_{m,0}$ note that the difference will be in having extra
canceling pairs of 1- and 2-handles that are again separated from the most of
other handles of $M_{m,0}$ and, as in the case of $Y_2$, can be removed from
the picture.

To complete the proof, we have to deal with the changes necessary to embed
more general $(U')^k_n$ into $M_{m,0}$.  Again we can use the trick from
Figures 10 and 11, by switching the roles of the biggest 0-framed and dotted
circles in the lover half of Figure 10 and by replacing the dot by 0-framing of
the two horizontal line segments throughout Figure 11.  Since all handle slides
were of 2-handles, the only difference is that we are sliding 2-handles over 1-
and 2-handles rather then 2-handles over only 1-handles and so we do not have
to use the forbidden link calculus moves involving slides of dotted circles
over components without a dot.\qed\enddemo

Note that although we have seen that any positive ribbon \rf\ can be associated
to a subsequence of the same sequence of stably non-product cobordisms $W_m$
from [BG], the particular $m$ needed to embed a given $U^k$ into $M_{m,1}$ and
$M_{m,0}$, and the number of their blow-ups, both depend on the given ribbon
\rf.

\enddocument